\documentclass[11pt,twoside,reqno]{amsart}
\usepackage{mathtools}
\usepackage{todonotes}
\usepackage{microtype}
\usepackage{cite}
\usepackage[OT1]{fontenc}
\usepackage{type1cm}
\usepackage{amssymb}
\usepackage{dsfont}
\usepackage{enumitem}
\usepackage{comment}
\usepackage{xcolor}
\usepackage{subfig}
\usepackage{aliascnt}
\usepackage{tikz}
\usetikzlibrary{calc,patterns,angles,quotes,arrows}
\usepackage{lmodern}

\usepackage{geometry}
\usepackage{hyperref}
\hypersetup{
  colorlinks=true,
  linkcolor=black,
  anchorcolor=black,
  citecolor=black,
  filecolor=black,      
  menucolor=red,
  runcolor=black,
  urlcolor=black,
}
\usepackage[capitalize,noabbrev,nameinlink]{cleveref}

\numberwithin{equation}{section}
\crefformat{equation}{(#2#1#3)}
\crefformat{figure}{#2Figure #1#3}
\crefformat{enumi}{(#2#1#3)}
\crefformat{theorem}{#2Theorem~#1#3}
\crefformat{lemma}{#2Lemma~#1#3}
\crefformat{proposition}{#2Proposition~#1#3}
\crefformat{corollary}{#2Corollary~#1#3}
\crefformat{claim}{#2Claim~#1#3}
\crefformat{section}{#2Section~#1#3}
\crefformat{subsection}{#2§#1#3}

\theoremstyle{plain}
\newtheorem{theorem}{Theorem}[section]

\newaliascnt{corollary}{theorem}

\aliascntresetthe{corollary}
\newaliascnt{proposition}{theorem}
\newtheorem{proposition}[proposition]{Proposition}
\aliascntresetthe{proposition}

\newaliascnt{lemma}{theorem}
\newtheorem{lemma}[lemma]{Lemma}
\aliascntresetthe{lemma}
\newaliascnt{conjecture}{theorem}

\aliascntresetthe{conjecture}

\theoremstyle{remark}

\theoremstyle{definition}

\newtheorem*{definition*}{Definition}

\newcommand{\ttt}{\mathbf{t}}

\newcommand{\LL}{\mathcal{L}}

\renewcommand{\SS}{\mathcal{S}}

\newcommand{\WW}{\mathcal{W}}

\newcommand{\R}{\mathbb{R}}

\newcommand{\Z}{\mathbb{Z}}
\newcommand{\N}{\mathbb{N}}

\newcommand{\fii}{\varphi}
\newcommand{\roo}{\varrho}

\newcommand{\dd}{\,\mathrm{d}}
\renewcommand{\ge}{\geqslant}
\renewcommand{\le}{\leqslant}
\renewcommand{\geq}{\geqslant}
\renewcommand{\leq}{\leqslant}

\DeclareMathOperator{\udimm}{\overline{dim}_M}

\DeclareMathOperator{\dimh}{dim_H}

\DeclareMathOperator{\dist}{dist}

\DeclareMathOperator{\proj}{proj}

\DeclareMathOperator{\sgn}{sgn}

\DeclareMathOperator{\id}{Id}
\DeclareMathOperator{\graph}{graph}

\begin{document}

\title{Level sets of prevalent Weierstrass functions}

\author{Zolt\'an Buczolich}
\address[Zolt\'an Buczolich]
        {Department of Analysis \\
         Elte E\"otv\"os Lor\'and University \\
         P\'azm\'any P\'eter S\'et\'any 1/C \\ 
         1117 Budapest
         Hungary }
\email{zoltan.buczolich@ttk.elte.hu}

\author{Antti K\"aenm\"aki}
\address[Antti K\"aenm\"aki]
        {University of Eastern Finland \\ 
         Department of Physics and Mathematics \\
         P.O.\ Box 111 \\ 
         FI-80101 Joensuu \\ 
         Finland}
\email{antti@kaenmaki.net}

\author{Bal\'azs Maga}
\address[Bal\'azs Maga]
        {HUN-REN Alfr\'ed R\'enyi Institute of Mathematics \\
         Budapest \\
         Hungary}
\email{magab@renyi.hu}
\thanks{B. Maga was supported by the KKP 139502 project, funded by the Ministry of Innovation and Technology of Hungary from the National Research, Development and Innovation Fund.}

% 03E15 - Descriptive set theory
% 26A16 - Lipschitz (H\"older) classes
% 26A21 - Classification of real functions; Baire classification of sets and functions
% 26B35 - Special properties of functions of several variables, H\"older conditions, etc
% 26E15 - Banach spaces of continuous, differentiable or analytic functions
% 28A50 - Integration and disintegration of measures
% 28A78 - Hausdorff and packing measures
% 54E52 - Baire category, Baire spaces
% 46E35 - Sobolev spaces and other spaces of ``smooth'' functions, embedding theorems, trace theorems
% 60G17 - Sample path properties
\date{\today}
%\thanks{}

\begin{abstract} 
  The $\alpha$-Weierstrass function is defined as $W_g^{\alpha,b}(x) = \sum_{k=0}^{\infty} b^{-\alpha k} g(b^k x)$, where $g$ is a Lipschitz function on the unit circle. For a prevalent $\alpha$-Weierstrass function, we prove that the upper Minkowski dimension of every level set is at most $1-\alpha$, and the Hausdorff dimension of almost every level set equals $1-\alpha$ with respect to its occupation measure. We further demonstrate that the occupation measure of a prevalent $\alpha$-Weierstrass function is absolutely continuous with respect to the Lebesgue measure. Consequently, the result on the Hausdorff dimension of level sets applies to a set of level sets with positive Lebesgue measure. A central tool in our analysis is the Weierstrass embedding. For a sufficiently large dimension $d$, we construct Lipschitz functions $g_0, g_1, \ldots, g_{d-1}$ such that the mapping $x \mapsto \big(W_{g_0}^{\alpha,b}(x), W_{g_1}^{\alpha,b}(x), \ldots, W_{g_{d-1}}^{\alpha,b}(x)\big)$ is $\alpha$-bi-H\"older. We also prove that such an embedding requires at least $1/\alpha$ coordinate functions.
\end{abstract}

\maketitle

\section{Introduction}\label{sec:intro}
A function $f\colon [0,1]\to \R$ is \emph{$\alpha$-H\"older} if there exists a constant $C>0$ such that
\begin{equation*}
  |f(x)-f(y)|\leq C|x-y|^\alpha
\end{equation*}
for all $x,y\in [0,1]$. The $\alpha$-H\"older condition implies that the upper Minkowski dimension of the graph of $f$ is at most $2-\alpha$; see \cite{MR1102677}. For the classical Weierstrass function, the H\"older condition is sharp, as the Minkowski dimension exists and equals $2-\alpha$; see, e.g., \cite[§11]{MR1102677}. While the Minkowski dimension of the graph is relatively straightforward to compute using elementary methods, determining the Hausdorff dimension usually requires more sophisticated techniques. The long-standing problem of computing the Hausdorff dimension of the Weierstrass function's graph was recently resolved by Shen \cite{Shen2018}. Additionally, Ren and Shen \cite{RenShen2021} established a dichotomy result for more general Weierstrass-type functions. These proofs rely heavily on the dynamics of the function graphs, utilizing careful transversality estimates in \cite{Shen2018} and powerful entropy methods introduced by Hochman \cite{Hochman2014} in \cite{RenShen2021}.

A subset $A$ of a Banach space $X$ is \emph{shy} if there exists a Borel measure $\mu$ and a Borel set $B \subset X$ containing $A$ such that $0 < \mu(K) < \infty$ for some compact set $K \subset X$ and $\mu(x+B) = 0$ for all $x \in X$. A set $A$ is \emph{prevalent} if its complement is shy. The notion of prevalence is abstract, but a concrete method to prove that a set $A \subset X$ is prevalent involves identifying a $d$-dimensional subspace $S \subset X$, called the \emph{probe space}, such that for any $x \in X$, we have $x + s \in B$ for Lebesgue-almost every $s \in S$, where the Lebesgue measure
$\LL^{d}$ is defined on $S$ and $B$ is a Borel set contained in $A$. In this case, $A$ is said to be \emph{$d$-prevalent}. It follows directly from the definitions that a set with a probe space is prevalent: indeed, taking $\mu$ as the Lebesgue measure on $S$ and $K$ as  the closed unit ball  in $S$, we have
\begin{equation*}
  \mu(x + (X \setminus B)) = \mu(\{s \in S : s - x \notin B\}) = 0
\end{equation*}
for all $x \in X$.

For an integer $b \ge 2$, $0<\alpha\leq 1$, and a non-constant $\Z$-periodic Lipschitz function $g$ on $\mathbb{S}^1$, where $\mathbb{S}^1$ is identified with $\R/\Z$, the \emph{Weierstrass function} $W_g^{\alpha,b} \colon \mathbb{S}^1 \to \R$ is defined by
\begin{equation}\label{*Wdef}
  W_g^{\alpha,b}(x) = \sum_{k=0}^{\infty} b^{-\alpha k} g(b^kx).
\end{equation}
For fixed $b \ge 2$ and $0 < \alpha \leq 1$, the family of such Weierstrass functions is denoted by $\WW^{\alpha,b}$. These functions are known to be $\alpha$-H\"older (see, e.g., \cite[Proposition~2.3]{10.1007/978-3-319-18660-3_5}), implying that
\begin{equation} \label{eq:udimm-upper-bound}
  \udimm(\graph(W_g^{\alpha,b})) \le 2-\alpha,
\end{equation}
where $\udimm(A)$ denotes the upper Minkowski dimension of a bounded set $A$. The functions in $\mathcal{W}^{\alpha,b}$ are in one-to-one correspondence with the Banach space $\mathrm{Lip}(\mathbb{S}^1)$ of Lipschitz functions on $\mathbb{S}^1$, equipped with the norm
\begin{equation*}
  \|g\|_{\mathrm{Lip}} = \sup_{x \in \mathbb{S}^1}|g(x)| + \sup_{(x,y) \in \mathbb{S}^1\times \mathbb{S}^1} \frac{|g(x)-g(y)|}{|x-y|}.
\end{equation*}
Thus, we identify $\WW^{\alpha,b}$ with $\mathrm{Lip}(\mathbb{S}^1)$, and prevalence in $\WW^{\alpha,b}$ refers to prevalence under this identification. The family $\WW^{\alpha,b}$ includes well-studied examples such as the classical Weierstass function, where $g(x) = \cos(2\pi x)$, and the Takagi function, where $g(x) = \dist(x,\Z)$.

The level sets of the Takagi function have been recently investigated in \cite{AnttilaBaranyKaenmaki2023}. Earlier results include \cite{BZirr}, with comprehensive surveys provided in \cite{AlKa} and \cite{Lags}; see also the references therein for further details. We investigate the size of level sets for prevalent Weierstrass functions, defined as the horizontal slices of the function's graph or, equivalently, the fibers of the projection of the graph onto the $y$-axis. We denote the Hausdorff dimension of a set $A \subset \R^d$ by $\dimh(A)$. In the absence of obvious obstructions, the Hausdorff dimension of the level set $(W_g^{\alpha,b})^{-1}(\{y\})$ is expected to be at most $\dimh(\graph(W_g^{\alpha,b}))-1$. Indeed, by Marstrand's slicing theorem (see, e.g., \cite[Theorem 1.6.1]{BishopPeres2016}) and \cref{eq:udimm-upper-bound}, we have
\begin{equation} \label{eq:marstrand}
\begin{split}
  \dimh((W_g^{\alpha,b})^{-1}(\{y\})) &= \dimh(\graph(W_g^{\alpha,b}) \cap \proj_{2}^{-1}(\{y\})) \\ 
  &\le \dimh(\graph(W_g^{\alpha,b}))-1 \le 1-\alpha
\end{split}
\end{equation}
for $\LL^1$-almost every $y \in \R$, where $\proj_{2} \colon \R^2 \to \R$, $\proj_{2}(x,y) = y$, is the orthogonal projection onto the $y$-axis. Furthermore, the main theorem of \cite{RenShen2021} determines this dimension when $g$ is real analytic. There the function is written as $\sum_{k=0}^{\infty}\lambda^{k}g(b^{k}x)$ with an integer $b \geq 2$ and $\lambda \in (1/b,1)$, which is \cref{*Wdef} with $\lambda = b^{-\alpha}$ and $\alpha \in (0,1)$, and the result asserts that for every $\Z$-periodic real analytic $g$ either $W_g^{\alpha,b}$ is real analytic or the Hausdorff dimension of its graph equals $2+\log_{b}\lambda = 2-\alpha$. Moreover, if $g$ is non-constant, then the former alternative occurs only for finitely many $\lambda$, and hence $\dimh(\graph(W_g^{\alpha,b})) = 2-\alpha$ for all but finitely many $\alpha \in (0,1)$, the exceptional values depending on both $g$ and $b$.

The level sets of prevalent $\alpha$-H\"older functions were recently investigated in \cite{AnttilaBaranyKaenmaki2025}. The main result, stated in \cite[Theorem~1.6]{AnttilaBaranyKaenmaki2025}, is the following:

\begin{theorem} \label{thm:main_holder}
  A prevalent $\alpha$-H\"older function $f$ on the unit interval satisfies
  \begin{enumerate}
    \item\label{it:main1_holder} $\udimm(f^{-1}(\{y\})) \le 1-\alpha$ for all $y \in \R$ provided that $0<\alpha<\tfrac12$,
    \item\label{it:main2_holder} $\LL^1(\{y \in f([0, 1]): \dimh(f^{-1}(\{y\})) = 1-\alpha\}) > 0$ provided that $0<\alpha<1$.
  \end{enumerate}
\end{theorem}

It was conjectured in \cite[Conjecture~1.7]{AnttilaBaranyKaenmaki2025} that prevalent $\alpha$-Weierstrass functions exhibit similar behavior and, specifically, that $\dimh((W_g^{\alpha,b})^{-1}(\{y\})) = 1-\alpha$ for $\LL^1$-almost all $y \in W_g^{\alpha,b}(\mathbb{S}^1)$. We partially confirm this conjecture by establishing an exact analogue of \cref{thm:main_holder} for prevalent $\alpha$-Weierstrass functions in the following two theorems:

\begin{theorem}\label{thm:main1}
  For any integer $b \ge 2$, a prevalent function $g \in \mathrm{Lip}(\mathbb{S}^1)$ satisfies
  \begin{equation*}
    \udimm((W_g^{\alpha,b})^{-1}(\{y\})) \le 1-\alpha
  \end{equation*}
  for all $y \in \R$ provided that $0<\alpha<\tfrac12.$
\end{theorem}

Let $\mathcal{G}=\{g_0,\ldots,g_{d-1}\}$ be a finite collection of Lipschitz functions on $\mathbb{S}^1$. The \emph{Weierstrass embedding} $\Phi_{\mathcal{G}}^{\alpha,b}\colon \mathbb{S}^1\to\R^d$ associated to the collection $\mathcal{G}$, integer $b \ge 2$, and $0<\alpha<1$ is defined by
\begin{equation}\label{*Web}
  \Phi_{\mathcal{G}}^{\alpha,b}(x)=(W_{g_0}^{\alpha,b}(x),W_{g_1}^{\alpha,b}(x),\ldots,W_{g_{d-1}}^{\alpha,b}(x)).
\end{equation}
The discussion in \cite[Section 4]{AnttilaBaranyKaenmaki2025} suggests that to resolve the conjecture for prevalent $\alpha$-Weierstrass functions, it is sufficient to construct an $\alpha$-bi-H\"older Weierstrass embedding and apply the results from \cite{AnttilaBaranyKaenmaki2025}. For \cref{thm:main1}, this approach is successful: in \cref{sec:embedding}, we establish the existence of such embedding, providing a positive resolution to \cite[Question 4.1]{AnttilaBaranyKaenmaki2025} and facilitating a direct proof of \cref{thm:main1}. 

However, the previous method does not extend to the second part of the conjecture. Specifically, in proving \cref{thm:main_holder}\cref{it:main2_holder} the idea was first to analyze slices in almost every direction and then adjust the slice direction by adding a linear function to the H\"older function. In the Banach space $\WW^{\alpha,b}$, which consists of $1$-periodic functions, this adjustment is not feasible. Consequently, in addition to constructing an $\alpha$-bi-H\"older Weierstrass embedding, we must directly analyze the level sets.

\begin{theorem}\label{thm:main2}
  For any integer $b \ge 2$, a prevalent function $g \in \mathrm{Lip}(\mathbb{S}^1)$ satisfies
  \begin{equation*}
    \LL^1(\{y \in W_g^{\alpha,b}(\mathbb{S}^1) : \dimh((W_g^{\alpha,b})^{-1}(\{y\})) = 1-\alpha\}) > 0
  \end{equation*}
	provided that $0<\alpha<1$.
\end{theorem}

Our method in \cref{sec:occupation} streamlines the proof of \cref{thm:main_holder}\cref{it:main2_holder}, showing that for prevalent $\alpha$-H\"older functions, the result can be established without relying on the ``almost every rotation'' argument. To prove our main result in \cref{sec:occupation}, we analyze the occupation measures of Weierstrass functions. Occupation measures are commonly studied in the context of stochastic processes. For one-dimensional Brownian motion $B$, the occupation measure, defined as $\mu(A) = \mathcal{L}^1(B^{-1}(A))$, is almost surely absolutely continuous with respect to the Lebesgue measure, satisfying the \emph{local time (LT)} condition; see, e.g., \cite[Theorem~3.26]{Morters_Peres_2010}.
Geman and Horowitz \cite{GeHo} noted: 
\begin{quote}
  ``$\ldots$results can be successfully applied to random functions and fields, it is difficult to apply them to particular nonrandom functions. For example, an interesting open problem is \emph{to determine which functions representable as Fourier series} (for instance) \emph{are (LT)}$\ldots$''
\end{quote}
Bertoin \cite{Berom,Berhdlev} analyzed the occupation measures and Hausdorff dimensions of level sets for certain self-affine functions. Depending on parameter values, these functions either satisfy the (LT) condition or have a singular occupation measure.

In \cite{BZirr}, it is shown that the occupation measure of the Takagi function is singular, and Lebesgue almost every level set is finite. The Weierstrass-Cellerier function, defined as $W(x) = \sum_{n=0}^{\infty} 2^{-n} \sin(2\pi 2^n x)$, is also examined. It is established that the occupation measure of $W(x) + cx$ is singular for almost every $c \in \mathbb{R}$. This result is strengthened in \cite{BZom}, which demonstrates that there are no exceptional $c$ values, including the case $c = 0$, confirming that the occupation measure of $W$ is singular. Notably, the Weierstrass-Cellerier function and the Takagi function belong to the class $\mathcal{W}^{1,2}$, whereas our \cref{*prabscb} addresses prevalent functions in the spaces $\mathcal{W}^{\alpha,b}$, with $0 < \alpha < 1$, showing that these functions typically have absolutely continuous occupation measures.

The article is organized as follows: In \cref{sec:embedding}, we establish the existence of an $\alpha$-bi-H\"older Weierstrass embedding, followed by a closer inspection of its properties in \cref{*propaaabound} and the proof of \cref{thm:main1}. In \cref{sec:occupation}, we first demonstrate in \cref{*prabscb} the absolute continuity of the occupation measure for prevalent Weierstrass functions. We then develop a method to directly analyze level sets, which we apply to prove \cref{thm:main2}. These techniques further enable us to prove \cref{*propaaabound}. The paper concludes with \cref{secq}, presenting open questions.

\section{Weierstrass embedding}\label{sec:embedding}

Recall that $\mathbb{S}^1$ is identified with $\R/\Z$ and, for $x,y \in \mathbb{S}^1$, write $|x-y|$ for the distance between $x$ and $y$ on $\mathbb{S}^1$. We say that $\Phi \colon \mathbb{S}^1 \to \R^d$ is an \emph{$\alpha$-bi-H\"older map} if there are constants $c_1,c_2>0$ such that
\begin{equation} \label{eq:bi-holder-def}
  c_1 \le \frac{\|\Phi(x)-\Phi(y)\|}{|x-y|^\alpha} \le c_2
\end{equation}
for all $x,y \in \mathbb{S}^1$ with $x \ne y$. For simplicity, we choose $\|\cdot\|$ to be the $\ell^\infty$ norm. The key thing to note is that if one can construct an $\alpha$-bi-H\"older embedding $\Phi\colon \mathbb{S}^1 \to \R^d$, whose coordinate functions are in $\WW^{\alpha,b}$, then \cref{thm:main1} holds for prevalent functions in $\WW^{\alpha,b}$.

Before stating \cref{prop:bi-holder} we need two lemmas which take care of the cases when $x$ and $y$ are sufficiently apart.

\begin{lemma} \label{lem:local-separation}
  Suppose $b \ge 2$ is an integer and $0 < \alpha < 1$. For every distinct $x,y \in \mathbb{S}^1$ there are $\delta_{x,y} > 0$ and a Lipschitz function $g$, which is therefore $\alpha$-H\"older, such that
  \begin{equation} \label{eq:local-separation}
    |W_g^{\alpha,b}(u)-W_g^{\alpha,b}(v)| \ge \frac{1-b^{-\alpha}}{2}
  \end{equation}
  for all $(u,v) \in B(x,\delta_{x,y})\times B(y,\delta_{x,y})$.
\end{lemma}

\begin{proof}
  Choose $N$ such that
  \begin{equation} \label{eq:local-separation-tail}
    \sum_{k=N}^{\infty}b^{-\alpha k} < \frac{1-b^{-\alpha}}{2}.
  \end{equation}
  Whenever $\delta_0$ is chosen below, we also require $\delta_0 < \tfrac14$ and define
  \begin{equation*}
    g(z) =
    \begin{cases}
      1, &\text{if } |z-x| \le \delta_0, \\
      2-\delta_0^{-1}|z-x|, &\text{if } \delta_0 < |z-x| < 2\delta_0, \\
      0, &\text{if } |z-x| \ge 2\delta_0.
    \end{cases}
  \end{equation*}
  The function $g \colon \mathbb{S}^1 \to [0,1]$ is nonconstant and Lipschitz. Set $T_b(z) = bz \bmod 1$. Thus, $T_b^0(z) = z$ and $T_b^{k+1}(z) = T_b(T_b^k(z))$.

  First suppose that $T_b^k(y) \ne x$ for every $k$. Since only finitely many iterates are involved, we may choose $\delta_0 > 0$ such that $T_b^k(B(y,\delta_0))\cap B(x,2\delta_0) = \emptyset$ for every $k \in \{0,\ldots,N-1\}$. Set $\delta_{x,y} = \delta_0$. For every $v \in B(y,\delta_{x,y})$, the definition of $W_g^{\alpha,b}$ in \cref{*Wdef} and \cref{eq:local-separation-tail} give
  \begin{equation} \label{eq:local-separation-y}
    W_g^{\alpha,b}(v) = \sum_{k=0}^{\infty}b^{-\alpha k}g(b^kv) = \sum_{k=N}^{\infty}b^{-\alpha k}g(b^kv) \le \sum_{k=N}^{\infty}b^{-\alpha k} < \frac{1-b^{-\alpha}}{2}.
  \end{equation}
  For every $u \in B(x,\delta_{x,y})$, nonnegativity of $g$ gives
  \begin{equation} \label{eq:local-separation-x}
    W_g^{\alpha,b}(u) \ge g(u) = 1.
  \end{equation}
  Combining \cref{eq:local-separation-y} and \cref{eq:local-separation-x}, we obtain
  \begin{equation*}
    |W_g^{\alpha,b}(u)-W_g^{\alpha,b}(v)| \ge 1-\frac{1-b^{-\alpha}}{2} = \frac{1+b^{-\alpha}}{2} > \frac{1-b^{-\alpha}}{2},
  \end{equation*}
  which proves \cref{eq:local-separation} in this case.

  Now suppose that the forward orbit of $y$ contains $x$, and let $k_0 \ge 1$ be the first hitting time. Then
  \begin{equation} \label{eq:first-orbit-hit}
    T_b^{k_0}(y) = x,
  \end{equation}
  while $T_b^k(y) \ne x$ for every $k \in \{0,\ldots,k_0-1\}$. Choose $\delta_0 > 0$ so that $\{T_b^k(y) : k \in \{0,\ldots,k_0-1\}\}\cap B(x,2\delta_0) = \emptyset$. It follows that
  \begin{equation*}
    W_g^{\alpha,b}(y) = \sum_{k=0}^{k_0-1}b^{-\alpha k}g(b^ky)+b^{-\alpha k_0}W_g^{\alpha,b}(x) = b^{-\alpha k_0}W_g^{\alpha,b}(x).
  \end{equation*}
  Since $W_g^{\alpha,b}(x) \ge 1$, we have
  \begin{equation*}
    W_g^{\alpha,b}(x)-W_g^{\alpha,b}(y) \ge 1-b^{-\alpha k_0} \ge 1-b^{-\alpha}.
  \end{equation*}
  The last lower bound is twice the right-hand side of \cref{eq:local-separation}. Therefore, continuity of $W_g^{\alpha,b}$ gives $\delta_{x,y} > 0$ such that \cref{eq:local-separation} holds for all $(u,v) \in B(x,\delta_{x,y})\times B(y,\delta_{x,y})$.
\end{proof}

The function provided by \cref{lem:local-separation} depends on the pair of points it separates. The following lemma removes this dependence at the cost of passing to a finite collection: a single collection of Lipschitz functions separates every pair of points lying at least a fixed distance apart, and only the choice of index within the collection is left to depend on the pair.

\begin{lemma}\label{*lemD}
  Suppose $b \ge 2$ is an integer, $0<\alpha<1$ and $D>0$ is given. Then there exist finitely many Lipschitz functions $g_{0},\ldots,g_{d-1}$ on $\mathbb{S}^1$, which are therefore $\alpha$-H\"older, such that for all $x,y\in \mathbb{S}^1$ with $|x-y|\geq D$ there exists $i\in \{0,\ldots,d-1\}$ such that $|W_{g_{i}}^{\alpha,b}(x)-W_{g_{i}}^{\alpha,b}(y)|\geq \frac{1-b^{-\alpha}}{2}.$
\end{lemma}

\begin{proof} 
  By compactness of $C = \{(x,y) \in \mathbb{S}^1 \times \mathbb{S}^1 : |x-y| \ge D\}$ and \cref{lem:local-separation}, we can select finitely many functions $g_i$, $i \in \{0,\ldots,d-1\}$, such that \cref{eq:local-separation} holds for any $(x,y) \in C$ with a suitable $i$ depending on $(x,y)$.
\end{proof}

The following theorem is the main result of this section, and it answers \cite[Question 4.1]{AnttilaBaranyKaenmaki2025} in the affirmative. Since each coordinate function is $\alpha$-H\"older, the upper bound in \cref{eq:bi-holder-def} holds automatically and the difficulty lies entirely in the lower bound. \cref{*lemD} disposes of pairs of points that are a definite distance apart, so the proof concentrates on nearby pairs, for which we construct an explicit family of translated piecewise linear functions.

\begin{theorem}\label{prop:bi-holder}
  For every integer $b \ge 2$ and $0<\alpha<1$ there exist $d \in \N$ and a finite collection $\mathcal{G}=\{g_0,\ldots,g_{d-1}\}$ of Lipschitz functions on $\mathbb{S}^1$ such that the Weierstrass embedding $\Phi_{\mathcal{G}}^{\alpha,b}\colon \mathbb{S}^1\to\R^d$ is $\alpha$-bi-H\"older.
\end{theorem}

\begin{proof}
  Let $0<\ell\le\tfrac12$ and define a function $\overline{h}_\ell \colon \R \to \R$ by setting $\overline{h}_\ell(x) = h_\ell(x-[x])$ for all $x \in \R$, where $[x]$ is the integer part of $x$ and $h_\ell \colon [0,1) \to \R$ satisfies
  \begin{equation*}
    h_\ell(x) =
    \begin{cases}
      -\frac{1}{1-\ell}x, &\text{if } 0 \le x < \frac12(1-\ell), \\ 
      \frac{1}{\ell}x - \frac{1}{2\ell}, &\text{if } \frac12(1-\ell) \le x < \frac12(1+\ell), \\ 
      -\frac{1}{1-\ell}x + \frac{1}{1-\ell}, &\text{if } \frac12(1+\ell) \le x < 1
    \end{cases}
  \end{equation*}
  for all $x \in [0,1)$. Notice that $\overline{h}_\ell$ is a piecewise linear function with $-\tfrac12 \le \overline{h}_\ell(x) \le \tfrac12$ for all $x \in \R$. In the open sets $(\tfrac12(1-\ell)+k, \tfrac12(1+\ell)+k)$, $k \in \Z$, where $\overline{h}_\ell$ is increasing, we have $\overline{h}_\ell'(x) = \ell^{-1} \ge 2$, and in the open sets $(\tfrac12(1+\ell)+k, \tfrac12(1-\ell)+k+1)$, $k \in \Z$, where $\overline{h}_\ell$ is decreasing, we have 
  \begin{equation}\label{*hlest}
 -2 \le \overline{h}_\ell'(x) = -(1-\ell)^{-1} < -1; 
 \end{equation}
  see \cref{fig:function-h} for an illustration.

  \begin{figure}[t]
  \begin{tikzpicture}[scale=1.0]
    \begin{scope}
      \draw[->] (-0.5,0) -- (4.5,0) node [below=0.15] {\footnotesize $x$};
      \draw[->] (0,-1.5) -- (0,1.5) node [left=0.1] {\footnotesize $y$};
      \draw (0,0) -- (28/15,-1) -- (32/15,1) -- (4,0);
      \draw[dotted] (28/15,0) -- (28/15,-1);
      \draw[dotted] (32/15,0) -- (32/15,1);
      \draw[dotted] (0,1) -- (32/15,1);
      \draw[dotted] (0,-1) -- (28/15,-1);
      \filldraw[black] (28/15,0) circle (1pt);
      \filldraw[black] (32/15,0) circle (1pt);
      \filldraw[black] (4,0) circle (1pt);
      \filldraw[black] (0,1) circle (1pt);
      \filldraw[black] (0,-1) circle (1pt);
      \node at (1.35,0.35) {\footnotesize $\tfrac12(1-\ell)$};
      \node at (2.65,-0.35) {\footnotesize $\tfrac12(1+\ell)$};
      \node at (4,0.35) {\footnotesize $1$};
      \node at (-0.24,1) {\footnotesize $\tfrac12$};
      \node at (-0.35,-1) {\footnotesize $-\tfrac12$};
    \end{scope}
    \begin{scope}[shift={(7,0)}]
      \draw[->] (-0.5,0) -- (4.5,0) node [below=0.15] {\footnotesize $x$};
      \draw[->] (0,-1.5) -- (0,1.5) node [left=0.1] {\footnotesize $y$};
      \draw (0,-1/2) -- (14/15,-1) -- (18/15,1) -- (4,-1/2);
      \draw[dotted] (14/15,0) -- (14/15,-1);
      \draw[dotted] (18/15,0) -- (18/15,1);
      \draw[dotted] (0,1) -- (18/15,1);
      \draw[dotted] (0,-1) -- (14/15,-1);
      \draw[dotted] (4,0) -- (4,-1/2);
      \filldraw[black] (14/15,0) circle (1pt);
      \filldraw[black] (18/15,0) circle (1pt);
      \filldraw[black] (4,0) circle (1pt);
      \filldraw[black] (0,1) circle (1pt);
      \filldraw[black] (0,-1) circle (1pt);
      \node at (1.35-8/15,0.35) {\footnotesize $r_i$};
      \node at (2.65-19/15,-0.35) {\footnotesize $s_i$};
      \node at (4,0.35) {\footnotesize $1$};
      \node at (-0.24,1) {\footnotesize $\tfrac12$};
      \node at (-0.35,-1) {\footnotesize $-\tfrac12$};
    \end{scope}
  \end{tikzpicture}
  \caption{The function $h_\ell$ with $\ell=\tfrac{1}{15}$ is depicted on the left and the function $g_i$, where we have chosen $r_i=\tfrac{7}{30}$ and $s_i=\tfrac{9}{30}$ for illustrative purposes, is on the right.}
  \label{fig:function-h}
  \end{figure}
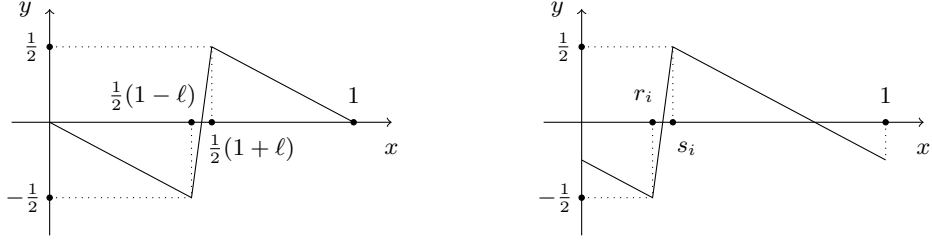
  Fix an integer $b \ge 2$ and $0<\alpha<1$. Let
  \begin{equation} \label{eq:definition-calpha}
    C_\alpha = \tfrac{10}{3} \alpha^{-1}\log_2 10 - \alpha^{-1}\log_2(1-2^{-\alpha}) + 2\alpha^{-1} - 1 > 1
  \end{equation}
  and choose $M \in \N$ such that
  \begin{equation} \label{eq:l0-choice}
    b^{M} > \max\biggl\{ \frac{20}{b^{1-\alpha}-1}, 40 C_\alpha b^{2C_\alpha(1-\alpha)} \biggr\}.
  \end{equation}
  Write $r_i = ib^{-M-3}$, $s_i = r_i+b^{-M}$, and $\ell = s_i-r_i = b^{-M} < \tfrac12$, and define a function $g_i \colon \mathbb{S}^1 \to \R$ by setting
  \begin{equation}\label{*gidef}
    g_i(x) = \overline{h}_{\ell}(x-\tfrac12(r_i+s_i-1))
  \end{equation}
  for all $x \in \mathbb{S}^1$ and $i \in \{0,\ldots,b^{M+3}-1\}$. Since $\overline{h}_\ell$ is $\Z$-periodic by construction, the right-hand side is independent of the representative chosen for $x$, so $g_i$ is well defined. Notice that each $g_i$ is a piecewise linear function with $-\tfrac12 \le g_i(x) \le \tfrac12$ for all $x \in \mathbb{S}^1$. In the open set $(r_i,s_i)$ (or $(0,s_i-1) \cup (r_i,1)$ if $s_i > 1$), where $g_i$ is increasing, we have $g_i'(x) = b^{M} \ge 2$, and when $g_i$ is decreasing, we have $-2 \le g_i'(x) = -(1-b^{-M})^{-1} < -1$.

  Let $d' = b^{M+3}$ and $\mathcal{G'} = \{g_i\}_{i=0}^{d'-1}$. Using  \cref{*lemD} with $D=b^{-M-1}$ we select finitely many functions  $\mathcal{G''} = \{g_i''\}_{i=0}^{d''-1}$ such that the conclusion of \cref{*lemD} holds for these functions. We put $\mathcal{G}=\mathcal{G'}\cup \mathcal{G''}$. Then we can write $\mathcal{G}= \{g_i\}_{i=0}^{d-1}$ with a suitable $d\geq d'$.
 
  It remains to show that the Weierstrass embedding $\Phi_{\mathcal{G}}^{\alpha,b}\colon \mathbb{S}^1\to\R^d$ defined in \cref{*Web} is $\alpha$-bi-H\"older. Since each $g_i$ is Lipschitz, each $W_{g_i}^{\alpha,b}$ is $\alpha$-H\"older, and hence the upper bound in \cref{eq:bi-holder-def} holds. For the lower bound, observe first that if $|x-y| \ge b^{-M-1}$, then the choice of $\mathcal{G''}$ gives an $i$ with $|W_{g_i}^{\alpha,b}(x)-W_{g_i}^{\alpha,b}(y)| \ge \tfrac{1-b^{-\alpha}}{2} \ge \tfrac{1-b^{-\alpha}}{2}|x-y|^{\alpha}$, where the last inequality uses $|x-y| \le 1$. Regarding each $W_{g_i}^{\alpha,b}$ as a $\Z$-periodic function on $\R$, we may choose the representatives of the remaining pairs so that $x<y$ and $y-x = |x-y|$. It therefore suffices to show that for each $x<y$ with $y-x < b^{-M-1}$ there is $i$ such that
  \begin{equation} \label{eq:embedding-goal}
    |W_{g_i}^{\alpha,b}(y)-W_{g_i}^{\alpha,b}(x)| \ge c|x-y|^\alpha,
  \end{equation}
  where the constant $c>0$ depends only on $b$ and $\alpha$. This is illustrated in \cref{fig:embedding}: We selected the function $g$ depicted in the upper left part of the figure and chose $b = 2$ and $\alpha = 0.7$ to define $W_g^{\alpha,b}$, shown in blue in the lower left part of the figure. We generated three translated copies of $g$, denoted 
  by
  $g_i$ for $i \in \{1, 2, 3\}$, with translations by $0$ (blue function), $0.214$ (red function), and $0.534$ (green function). The latter two values were chosen randomly to avoid unintended symmetry or correlation in the resulting images.

  The zero loci of $W_{g_i}^{\alpha,b}(y) - W_{g_i}^{\alpha,b}(x)$ for $i \in \{1, 2, 3\}$ are plotted in blue, red, and green, respectively. Regions violating \cref{eq:embedding-goal} are indicated with lighter shading. For a fixed $i$, the function $W_{g_i}^{\alpha,b}(y) - W_{g_i}^{\alpha,b}(x)$ equals zero on the diagonal $y = x$, so the diagonal is included in the zero level sets of all three functions, appearing blackish or grayish in the figure due to the overlap of all three colors. Since the right-hand side of \cref{eq:embedding-goal} is also zero when $x = y$, the presence of all three colors on the diagonal poses no issue.

  To satisfy \cref{eq:embedding-goal}, there must exist a constant $c$ such that no point off the diagonal is simultaneously shaded by all three colors. In the figure, we used $c = 0.2$. As $c \to 0$, the lightly shaded regions contract toward the zero loci, which can become highly complex near the diagonal. Given the resolution limits of our images, the absence of off-diagonal points shaded by all three colors does not conclusively confirm that a specific $c$ is suitable for these functions.

  If $W_{g_i}^{\alpha,b}(y) - W_{g_i}^{\alpha,b}(x) = 0$ at a pixel corresponding to $(x, y)$ off the diagonal for some $i$, then \cref{eq:embedding-goal} cannot be satisfied for that $i$ with any $c$. This represents a significant obstacle to achieving the Weierstrass embedding with these three functions.
  \begin{figure}%
  \centering
  \begin{minipage}{0.35\textwidth}%
        \includegraphics[width=\textwidth]{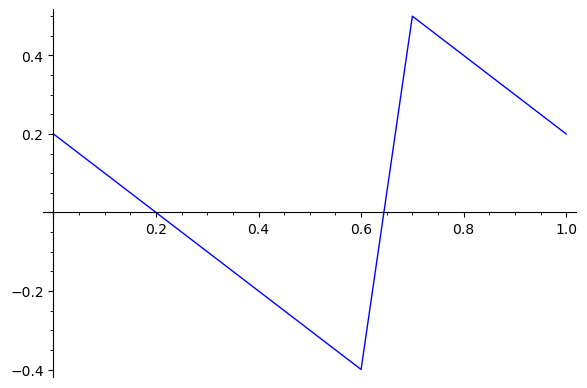}\\[20pt]
        \includegraphics[width=\textwidth]{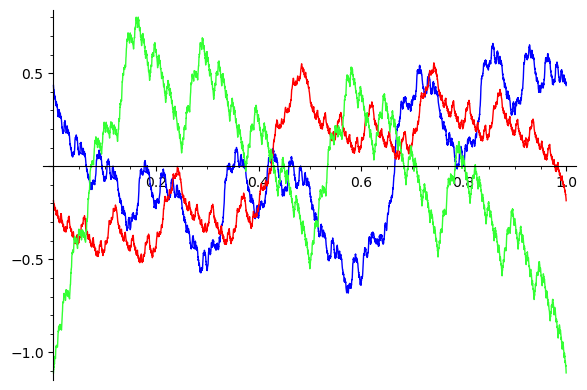}%
  \end{minipage}%
  \qquad
  \begin{minipage}{0.55\textwidth}%
        \includegraphics[width=\textwidth]{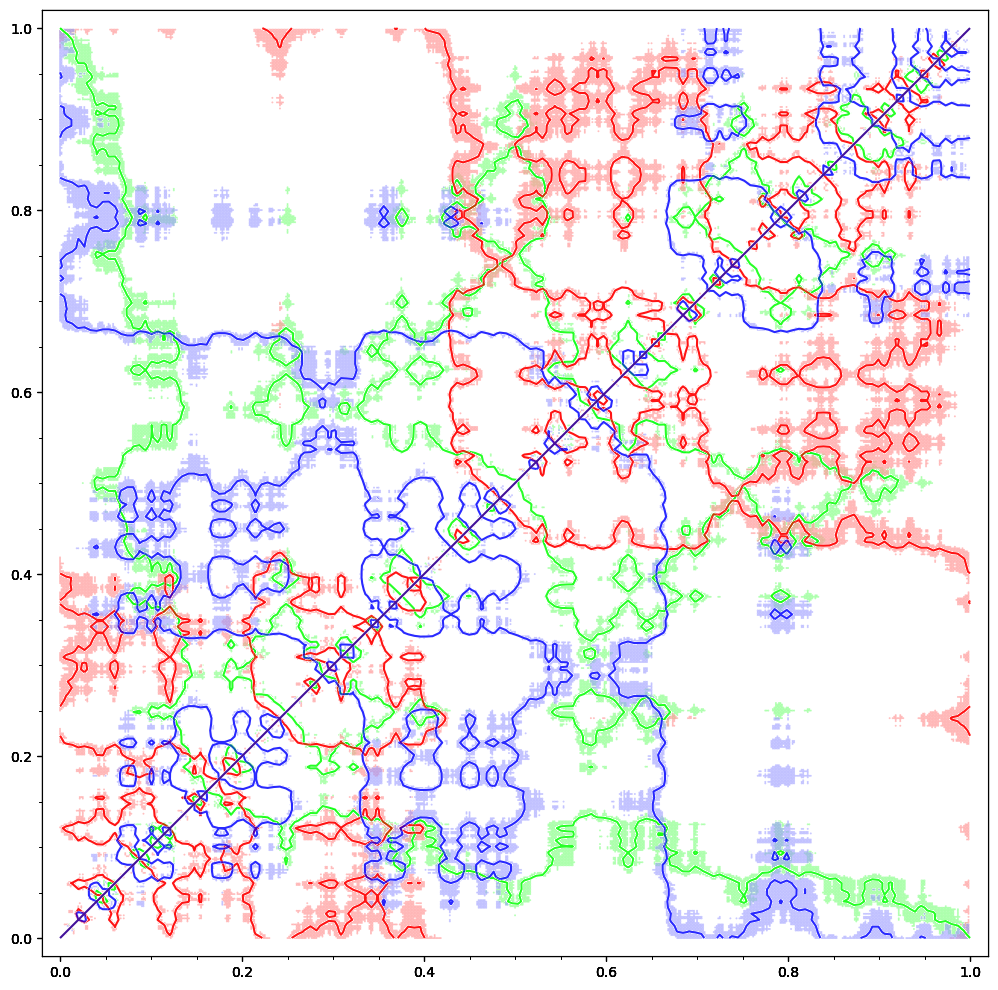}%
  \end{minipage}%
  \caption{The blue function in the bottom-left figure represents the Weierstrass function $W_g^{\alpha,b}$, corresponding to the function $g$ shown in the top-left figure. The right-hand figure depicts the zero loci and regions violating \cref{eq:embedding-goal} for three functions defined by $(x,y) \mapsto W_{g_i}^{\alpha,b}(y) - W_{g_i}^{\alpha,b}(x)$, where $i \in \{1, 2, 3\}$, and $g_1, g_2, g_3$ are distinct translations of $g$. All three functions are displayed in the bottom-left figure.}
  \label{fig:embedding}%
  \end{figure}
  To prove \cref{eq:embedding-goal} the idea is to find a dominant summand in the series defining the difference $W_{g_i}^{\alpha,b}(y)-W_{g_i}^{\alpha,b}(x)$ and then to show that the other summands do not contribute too much. Let $x<y$ be such that $y-x < b^{-M-1}$ and choose $k_{xy} \in \N$ 
  such that
  
  \begin{equation} \label{eq:kxy-choice}
    \frac{1}{b^{k_{xy}+M+2}} \le y-x < \frac{1}{b^{k_{xy}+M+1}}.
  \end{equation} 
  Let $j \in \N$ be such that
  \begin{equation*}
    \frac{j-1}{b^{k_{xy}}} < x \le \frac{j}{b^{k_{xy}}}
  \end{equation*}
  and choose $i \in \{0,\ldots,d-1\}$ such that the point $r_i = ib^{-M-3}$ satisfies
  \begin{equation} \label{eq:i-choice}
    \frac{j-1}{b^{k_{xy}}} + \frac{r_i}{b^{k_{xy}}} = \frac{j-1}{b^{k_{xy}}} + \frac{i}{b^{k_{xy}+M+3}} < x \le \frac{j-1}{b^{k_{xy}}} + \frac{i+1}{b^{k_{xy}+M+3}} = \frac{j-1}{b^{k_{xy}}} + \frac{r_{i+1}}{b^{k_{xy}}}.
  \end{equation}
  By the choices \cref{eq:kxy-choice} and \cref{eq:i-choice}, we see that
  \begin{align*}
    y &< x + \frac{1}{b^{k_{xy}+M+1}} \le \frac{j-1}{b^{k_{xy}}} + \frac{r_i}{b^{k_{xy}}} + \frac{1}{b^{k_{xy}+M+3}} + \frac{1}{b^{k_{xy}+M+1}} \\ 
    &< \frac{j-1}{b^{k_{xy}}} + \frac{r_i}{b^{k_{xy}}} + \frac{2}{b} \cdot \frac{1}{b^{k_{xy}+M}} \le \frac{j-1}{b^{k_{xy}}} + \frac{s_i}{b^{k_{xy}}},
  \end{align*}
  recalling that $s_i = r_i + b^{-M}$. The points $b^{k_{xy}}x \bmod 1$ and $b^{k_{xy}}y \bmod 1$ are thus both contained in the interval $[r_i,s_i]$ (or the union $[0,s_i-1] \cup [r_i,1)$ if $s_i \ge 1$). Since $g_i'(x) = b^{M}$ on the corresponding open interval, we have
  \begin{equation} \label{eq:difference-series1}
    b^{-\alpha k_{xy}} g_i(b^{k_{xy}}y) - b^{-\alpha k_{xy}} g_i(b^{k_{xy}}x) = b^{(1-\alpha)k_{xy}} b^{M}(y-x)
  \end{equation}
  The summand corresponding to the index $k_{xy}$ in the series defining the difference $W_{g_i}^{\alpha,b}(y)-W_{g_i}^{\alpha,b}(x)$ will give us the desired behavior. 

  Let us next show that the contribution of other summands (see \cref{*Wdef}) remains so small that they do not have an effect on the behavior of the series defining the difference $W_{g_i}^{\alpha,b}(y)-W_{g_i}^{\alpha,b}(x)$. Notice that if $k$ does not satisfy $g_i(b^ky)-g_i(b^kx) \ge 0$, then by the definition of $g_i$, namely by \cref{*hlest} and \cref{*gidef}, it satisfies
  \begin{equation} \label{eq:decreasing-k}
    g_i(b^ky)-g_i(b^kx) \ge -2b^k(y-x).
  \end{equation}
  By \cref{eq:l0-choice}, we have
  \begin{equation*}
    \sum_{k=0}^{k_{xy}-1} b^{(1-\alpha)k} = \frac{b^{(1-\alpha)k_{xy}}-1}{b^{1-\alpha}-1} < \tfrac{1}{20} b^{(1-\alpha)k_{xy}+M}.
  \end{equation*}
  Therefore,
  \begin{equation} \label{eq:difference-series2}
  \begin{split}
    \sum_{k=0}^{k_{xy}-1} b^{-\alpha k} g_i(b^ky) - b^{-\alpha k} g_i(b^kx) &\ge -2\sum_{k=0}^{k_{xy}-1} b^{(1-\alpha)k}(y-x) \\ 
    &> -\tfrac{1}{10} b^{(1-\alpha)k_{xy}}b^{M} (y-x).
  \end{split}
  \end{equation}
  To estimate the tail, recall the definition of $C_\alpha$ from \cref{eq:definition-calpha} and fix an integer $k_{xy}' \ge k_{xy}$ such that $C_\alpha < k_{xy}' - k_{xy} < 2C_\alpha$. Since $\log b \ge \log 2 > \tfrac{3}{10}$ and $-\log_b(1-b^{-\alpha})$ is decreasing in $b$, we have
  \begin{equation*}
    \alpha(k_{xy}' - k_{xy}) > \log_b 10 - \log_b(1-b^{-\alpha}) + 2 - \alpha
  \end{equation*}
  and hence,
  \begin{equation} \label{eq:kxy2-choice}
    \frac{b^{-\alpha(k_{xy}'-1)}}{1-b^{-\alpha}} < \tfrac{1}{10}b^{-\alpha k_{xy}-2}.
  \end{equation}
  Since $|g_i| \le \tfrac12$, the inequalities \cref{eq:kxy2-choice} and \cref{eq:kxy-choice} give
  \begin{equation} \label{eq:difference-series3}
  \begin{split}
    \biggl| \sum_{k=k_{xy}'+1}^\infty b^{-\alpha k}g_i(b^ky) - b^{-\alpha k}g_i(b^kx) \biggr| &\le \biggl| \sum_{k=k_{xy}'+1}^\infty b^{-\alpha k} \biggr|  
    < \tfrac{1}{10}b^{-\alpha k_{xy}-2} \\ &\le \tfrac{1}{10}b^{(1-\alpha)k_{xy}} b^{M} (y-x).
  \end{split}
  \end{equation}
  Furthermore, by the choice of $k_{xy}'$ and \cref{eq:l0-choice}, we have $40C_\alpha b^{(1-\alpha)(k_{xy}'-k_{xy})} \le b^{M}$ and
  \begin{equation}\label{*Caaa}
    4C_\alpha b^{(1-\alpha)k_{xy}'} < \tfrac{1}{10}b^{(1-\alpha)k_{xy}} b^{M}.
  \end{equation}
  Therefore, by \cref{eq:decreasing-k} and \cref{*Caaa}, we see that
  \begin{equation} \label{eq:difference-series4}
  \begin{split}
    \sum_{k=k_{xy}+1}^{k_{xy}'} b^{-\alpha k}g_i(b^ky) - b^{-\alpha k}g_i(b^kx) &> -2\sum_{k=k_{xy}+1}^{k_{xy}'} b^{(1-\alpha)k} (y-x) \\ 
    &\ge -2(k_{xy}'-k_{xy}) b^{(1-\alpha)k_{xy}'} (y-x) \\ 
    &\ge -4C_\alpha b^{(1-\alpha)k_{xy}'} (y-x) \\ 
    &> -\tfrac{1}{10}b^{(1-\alpha)k_{xy}} b^{M} (y-x).
  \end{split}
  \end{equation}
  Adding \cref{eq:difference-series1}, \cref{eq:difference-series2}, \cref{eq:difference-series3}, and \cref{eq:difference-series4} together and applying \cref{eq:kxy-choice}, we conclude that
  \begin{align*}
    W_{g_i}^{\alpha,b}(y)-W_{g_i}^{\alpha,b}(x) &\ge \tfrac{7}{10}b^{(1-\alpha)k_{xy}} b^{M} (y-x) \\ 
    &= \tfrac{7}{10}b^{-(1-\alpha)(-k_{xy}-M-2)} b^{-(1-\alpha)(M+2)} b^{M} (y-x) \\ 
    &\ge \tfrac{7}{10}b^{\alpha M-2(1-\alpha)} (y-x)^\alpha
  \end{align*}
  as claimed in \cref{eq:embedding-goal} with $c = \tfrac{7}{10}b^{\alpha M-2(1-\alpha)}$. Together with the case $|x-y| \ge b^{-M-1}$ treated above, this proves the lower bound in \cref{eq:bi-holder-def} with $c_1 = \min\bigl\{c,\tfrac{1-b^{-\alpha}}{2}\bigr\}$.
\end{proof}

The following proposition demonstrates that in \cref{prop:bi-holder}, the minimal number of Lipschitz functions on $\mathbb{S}^1$ required for the Weierstrass embedding $\Phi_{\mathcal{G}}^{\alpha,b} \colon \mathbb{S}^1 \to \R^d$ to be $\alpha$-bi-H\"older is at least $\frac{1}{\alpha}$. 

\begin{proposition} \label{*propaaabound}
  Let $b \geq 2$ and $d \geq 1$ be integers, and let $0 < \alpha < 1$. If $d < \frac{1}{\alpha}$ and $\mathcal{G} = \{g_0, \ldots, g_{d-1}\}$ is a finite collection of Lipschitz functions on $\mathbb{S}^1$, then the Weierstrass embedding $\Phi_{\mathcal{G}}^{\alpha,b} \colon \mathbb{S}^1 \to \mathbb{R}^d$ is not $\alpha$-bi-H\"older.
 \end{proposition}
 
 \begin{proof}
  Suppose for contradiction that $\Phi_{\mathcal{G}}^{\alpha,b}$ is $\alpha$-bi-H\"older. Then there exists a constant $c_1 > 0$ such that $\|\Phi_{\mathcal{G}}^{\alpha,b}(x)-\Phi_{\mathcal{G}}^{\alpha,b}(y)\| \ge c_1|x-y|^\alpha$ for all $x,y \in \mathbb{S}^1$. This makes $\Phi_{\mathcal{G}}^{\alpha,b}$ injective and its inverse $(\Phi_{\mathcal{G}}^{\alpha,b})^{-1} \colon \Phi_{\mathcal{G}}^{\alpha,b}(\mathbb{S}^1) \to \mathbb{S}^1$ satisfies $|(\Phi_{\mathcal{G}}^{\alpha,b})^{-1}(u)-(\Phi_{\mathcal{G}}^{\alpha,b})^{-1}(v)| \le c_1^{-1/\alpha}\|u-v\|^{1/\alpha}$, so $(\Phi_{\mathcal{G}}^{\alpha,b})^{-1}$ is $\frac1\alpha$-Hölder. A $\gamma$-Hölder map $F$ satisfies $\dimh(F(A)) \le \frac1\gamma\dimh(A)$ for every set $A$, by \cite[Proposition~2.3]{MR1102677}, so
  \begin{equation*}
    1 = \dimh(\mathbb{S}^1) = \dimh((\Phi_{\mathcal{G}}^{\alpha,b})^{-1}(\Phi_{\mathcal{G}}^{\alpha,b}(\mathbb{S}^1))) \le \alpha\dimh(\Phi_{\mathcal{G}}^{\alpha,b}(\mathbb{S}^1)) \le \alpha d ,
  \end{equation*}
  which is a contradiction since $d < \frac{1}{\alpha}$.
 \end{proof}

Let $\Phi = \Phi_{\mathcal{G}}^{\alpha,b}$ be the Weierstrass embedding associated with $\mathcal{G} = \{g_0,\ldots,g_{d-1}\}$, where $b\geq 2$ is an integer and $0<\alpha<1$. For any $W \in \WW^{\alpha,b}$ and $\ttt \in \R^d$, the function $W_{\ttt} \colon \mathbb{S}^1 \to \R$, defined by
\begin{equation*}
  W_{\ttt}(x)=W(x)+\langle\ttt,\Phi(x)\rangle,
\end{equation*}
belongs to $\WW^{\alpha,b}$. The space
\begin{equation*}
  \SS = \{x\mapsto \langle\ttt, \Phi(x)\rangle : \ttt \in \R^d\}
\end{equation*}
is a $d$-dimensional subspace of $\WW^{\alpha,b}$ and serves as the probe space in our results. To establish that a property holds for $d$-prevalent functions, it is enough to show that for any $W\in\WW^{\alpha,b}$, the property holds for $W_{\ttt}$ for $\LL^d$-almost every $\ttt\in\R^d$.

\begin{proof}[Proof of \cref{thm:main1}]
  Let $b \geq 2$ an integer, $0<\alpha<1$, and let $\Phi=\Phi_{\mathcal{G}}^{\alpha,b}$ be the $\alpha$-bi-H\"older Weierstrass embedding given by \cref{prop:bi-holder}. Although the map $\Phi$ is fixed in \cite[Section~3]{AnttilaBaranyKaenmaki2025} to be an $\alpha$-bi-H\"older map on $[0,1]$, the proof of \cite[Proposition~3.2]{AnttilaBaranyKaenmaki2025} uses it only through the lower bound $\|\Phi(x)-\Phi(y)\| \geq c_1|x-y|^{\alpha}$ and the resulting finiteness of $\int_0^1\int_0^1|x-y|^{-\alpha s} \dd x \dd y$ for $\alpha s<1$, both of which remain valid when $|x-y|$ denotes the distance on $\mathbb{S}^1$. By \cite[Proposition~3.2]{AnttilaBaranyKaenmaki2025}, for every $W\in\WW^{\alpha,b}$ and for $\LL^d$-almost every $\ttt \in \R^d$, the $\alpha$-H\"older function $W_{\ttt} \in \mathcal{W}^{\alpha,b}$ defined by $W_{\ttt}(x)=W(x)+\langle\ttt,\Phi(x)\rangle$ satisfies
  \begin{equation*}
    \udimm(W_{\ttt}^{-1}(\{y\})) \le 1-\alpha
  \end{equation*}
  for all $y \in \R$.
\end{proof}

\section{Fourier analysis and occupation measures}\label{sec:occupation}

Let $b \ge 2$ be an integer and $0<\alpha<1$. Let $\Phi = \Phi_{\mathcal{G}}^{\alpha,b}$ be the $\alpha$-bi-H\"older Weierstrass embedding associated with $\mathcal{G} = \{g_0,\ldots,g_{d-1}\}$, as given by \cref{prop:bi-holder}. For any $W \in \WW^{\alpha,b}$ and $\ttt \in \R^d$, the function $W_{\ttt} \colon \mathbb{S}^1 \to \R$, defined by
\begin{equation*}
  W_{\ttt}(x)=W(x)+\langle\ttt,\Phi(x)\rangle,
\end{equation*}
belongs to $\WW^{\alpha,b}$. The proof of \cref{thm:main2} employs Fourier-analytic and potential-theoretic methods. These tools are applied to the natural measures $\mu_\ttt$ supported on the graphs of the functions $W_{\ttt}$, defined as
\begin{equation} \label{eq:lift-measure}
  \mu_{\ttt}=(\id,W_{\ttt})_*\LL^1.
\end{equation}
Thus, $\mu_{\ttt}$ is the lift of the Lebesgue measure from the unit interval onto the graph of $W_{\ttt}$. Recalling that $\proj_{2} \colon \R^2 \to \R$, $\proj_{2}(x,y) = y$, is the projection onto the $y$-axis, the \emph{occupation measure} $\lambda_\ttt$ associated with $W_{\ttt}$ is given by
\begin{equation*}
  \lambda_\ttt = (\proj_{2})_*\mu_\ttt.
\end{equation*}
Equivalently, the occupation measure is the pushforward of the Lebesgue measure, $\lambda_\ttt = (W_{\ttt})_* \LL^1$.

The Fourier transform of the occupation measure $\lambda_\ttt$ is defined as
\begin{equation*}
  \widehat\lambda_\ttt(\xi) = \int_{\R} e^{i\xi x} \dd\lambda_\ttt(x)
\end{equation*}
for all $\xi \in \R$. The Fourier transform of an integrable function $\psi \colon \R^d \to \R$ is given by
\begin{equation*}
  \widehat\psi(\zeta) = \int_{\R^d} \psi(\ttt) e^{i\langle \ttt,\zeta \rangle} \dd\ttt
\end{equation*}
for all $\zeta \in \R^d$, where $\langle \,\cdot\,,\cdot\, \rangle$ denotes the standard inner product on $\R^d$. Let $C^{\infty}_0(\R^d)$ denote the space of smooth functions with compact support on $\R^d$. By the Riemann-Lebesgue lemma and integration by parts, for any $\psi \in C^{\infty}_0(\R^d)$ and $s>0$, there exists a constant $C=C(\psi,s)>0$ such that
\begin{equation} \label{eq:riemann-lebesgue}
  |\widehat\psi(\zeta)|\leq C(1+\|\zeta\|)^{-s}
\end{equation}
for all $\zeta\in\R^d$.

The following proposition establishes that the occupation measure $\lambda_\ttt$ is absolutely continuous for almost every $\ttt$. By the Radon-Nikodym theorem (see, e.g., \cite[Theorem~2.12(2)]{Mattila1995}), the measure $\lambda_\ttt$ is absolutely continuous if and only if it has a Radon-Nikodym derivative $\mathrm{d}\lambda_\ttt/\mathrm{d}y \in L^1(\R)$. 

\begin{theorem}\label{*prabscb}
  Let $b \ge 2$ be an integer, $0<\alpha<1$, and $W \in \WW^{\alpha,b}$. Then for $\LL^d$-almost every $\ttt \in \R^d$ the occupation measure $\lambda_\ttt$ associated with $W_\ttt$ is absolutely continuous with density in $L^2(\R)$.
\end{theorem}

It follows directly that $\alpha$-Weierstrass functions satisfying the local time (LT) condition are $d$-prevalent in $\mathcal{W}^{\alpha,b}$. The proof of \cref{*prabscb} builds on \cite[Proposition 3.2]{AnttilaBaranyKaenmaki2025}, which applies only for $0 < \alpha < \frac{1}{2}$ and yields a stronger result not explicitly stated but derived within its proof: namely, the occupation measure is not only absolutely continuous but also has a bounded and continuous density. By combining these techniques with our methods in the present work, we establish the following theorem.

\begin{theorem}\label{*prabscbb}
  Let $b \ge 2$ be an integer, $0<\alpha<\frac12$, and $W \in \WW^{\alpha,b}$. Then for $\LL^d$-almost every $\ttt \in \R^d$ the occupation measure $\lambda_\ttt$ associated with $W_\ttt$ is absolutely continuous with bounded and continuous density.
\end{theorem}

The proof of \cref{*prabscbb} combines the arguments from \cref{*prabscb} and \cite[Proposition 3.2]{AnttilaBaranyKaenmaki2025}, and we omit the details for brevity. Consequently, for $0<\alpha<\frac12$, \cref{*prabscbb} implies that $\alpha$-Weierstrass functions with an absolutely continuous occupation measure possessing a bounded and continuous density are $d$-prevalent in $\mathcal{W}^{\alpha,b}$.

\begin{proof}[Proof of \cref{*prabscb}]
  Since $0 < \alpha < 1$, choose  $s > 1$ such that  $\alpha s < 1$. Let $\roo>0$ and let $\psi\in C^\infty_0(\R^d)$ be a smooth function satisfying $0\leq \psi\leq 1$ and $\psi(\ttt)=1$ for all $\ttt\in B(\mathbf{0},\roo)$. By Fubini's theorem and \cref{eq:riemann-lebesgue}, there exists a constant $C = C(\psi,s)>0$ such that
  \begin{align*}
    \biggl| \int\limits_{B(\mathbf{0},\roo)} \int  &|\widehat\lambda_{\ttt}(\xi)|^2 \dd\xi\dd\textbf{t} \biggr| \leq \biggl| \iint\int_0^1\int_0^1 \psi(\textbf{t})  e^{i\xi(W(x)-W(y)+\langle \textbf{t},\Phi(x)-\Phi(y) \rangle)} \dd x\dd y\dd\xi\dd\textbf{t} \biggr| \\ 
    &= \biggl| \int\int_0^1\int_0^1  e^{i\xi(W(x)-W(y))}\Big ( \int \psi(\textbf{t}) e^{i\langle \textbf{t},\xi(\Phi(x)-\Phi(y)) \rangle} \dd\textbf{t}\Big )\dd x\dd y\dd\xi \biggr| \\ 
    &= \biggl| \int\int_0^1\int_0^1  e^{i\xi(W(x)-W(y))} \widehat\psi(\xi(\Phi(x)-\Phi(y))) \dd x\dd y\dd\xi \biggr| \\ 
    &\le \int\int_0^1\int_0^1  |\widehat\psi(\xi(\Phi(x)-\Phi(y)))| \dd x\dd y\dd\xi \\ 
    &\le \int\int_0^1\int_0^1 \frac{C}{(1+|\xi|\|\Phi(x)-\Phi(y)\|)^s} \dd x\dd y\dd\xi.
  \end{align*}
  Assuming that the constant $c_1 \leq 1$ in the $\alpha$-bi-H\"older condition \cref{eq:bi-holder-def}, and using that $\Phi$ is $\alpha$-bi-H\"older with $|x - y| \leq 1$, we obtain
  \begin{align*}
    \int\int_0^1\int_0^1& \frac{C}{(1+|\xi|\|\Phi(x)-\Phi(y)\|)^s} \dd x\dd y\dd\xi\leq\int\int_0^1\int_0^1 \frac{C}{(1+c_1|\xi||x-y|^{\alpha})^s} \dd x\dd y\dd\xi \\
    &\leq \int\int_0^1\int_0^1 \frac{C}{\left(\frac{1}{|x-y|^{\alpha}}+c_1|\xi|\right)^s|x-y|^{\alpha s}} \dd x\dd y\dd\xi\\
    &\leq\int\frac{C}{c_1^s(1+|\xi|)^s}\dd\xi\int_0^1\int_0^1\frac{1}{|x-y|^{\alpha s}} \dd x\dd y<\infty
  \end{align*}
  since $\alpha s < 1$. Thus, for $\LL^d$-almost every $\ttt$, the Fourier transform $\widehat{\lambda}_\ttt$ is in $L^2(\R)$. Since $\lambda_\ttt$ is a compactly supported finite measure, \cite[Theorem~3.3]{MattilaFou} shows that $\lambda_\ttt$ is absolutely continuous with density $\mathrm{d}\lambda_\ttt/\mathrm{d}y \in L^2(\R)$. As $\lambda_\ttt$ is a probability measure, this density belongs to $L^1(\R)$ as well.
\end{proof}

By \cref{*prabscb}, for $\LL^d$-almost every $\ttt \in \R^d$, there exists a function $h_{\ttt} \in L^2(\R)$ such that 
\begin{equation*}
  \lambda_{\ttt}(A) = \int_A h_\ttt(y)\dd y
\end{equation*}
for all Borel sets $A \subset \R$. Applying Rohlin's disintegration theorem (see, e.g., \cite[Theorem 2.1]{Simmons2012}) to the projection $\proj_{2} \colon \R^2 \to \R$, defined by $\proj_{2}(x,y) = y$, and the measure $\mu_{\ttt}$, we obtain, for $\lambda_\ttt$-almost every $y$, a Borel measure $(\mu_{\textbf{t}})_{y}$ supported on $\graph(W_{\ttt}) \cap \proj_{2}^{-1}(\{y\})$. By \cite{Simmons2012}, we have
\begin{equation*}
  \mu_{\ttt}(A) = \iint \mathds{1}_A(a) \dd(\mu_{\ttt})_{y}(a) \dd\lambda_{\ttt}(y)
\end{equation*}
for all Borel sets $A \subset \R^2$, where the parameter $a$ is in $\R^2$. For each $y$, define a Borel measure by $(\nu_{\ttt})_{y} = h_{\ttt}(y)(\mu_{\ttt})_{y}$, and note that
\begin{equation} 
  \mu_{\ttt}(A) = \iint \mathds{1}_A(a) \dd(\nu_{\ttt})_{y}(a) \dd y
\end{equation}
for all Borel sets $A \subset \R^2$. This implies that for $\mu_{\ttt}$-integrable functions $f$, we have
\begin{equation}\label{*disint}
  \int f(a)\dd\mu_{\ttt}(a)=\iint f(a)\dd(\nu_{\ttt})_{y}(a) \dd y.
\end{equation}
Furthermore, by \cite[Theorem 2.2]{Simmons2012}, for all $g \in L^1(\R)$, we have
for $\lambda_{\ttt}$ almost every $y$
\begin{equation} \label{eq:grslbcl6}
  \int g(a) \dd(\mu_{\ttt})_{y}(a) = \lim_{r \downarrow 0} \frac{\int_{B^\top(y,r)} g(a) \dd\mu_{\ttt}(a)}{\lambda_{\ttt}(B(y,r))} = \lim_{r \downarrow 0} \frac{1}{2r} \frac{\int_{B^\top(y,r)} g(a) \dd\mu_{\ttt}(a)}{h_{\ttt}(y)},
\end{equation}
where $B^\top(y,r) = \proj_{2}^{-1}(B(y,r))$. Since the measures $(\mu_{\textbf{t}})_{y}$ are supported on $\graph(W_{\ttt}) \cap \proj_{2}^{-1}(\{y\})$ for given $y_{0}\in \R$ and $r>0$ we can deduce from \eqref{*disint} that  
\begin{equation}\label{eq:grslbcl5}
\begin{split}
  \int_{B^\top(y_{0},r) } f(a) \dd\mu_{\ttt}(a) &= \int \mathds{1}_{B(y_{0},r)}(\proj_{2}(a)) f(a) \dd\mu_{\ttt}(a) \\ 
  &= \int_{y_{0}-r}^{y_{0}+r}\int f(a)\dd(\nu_{\ttt})_{y}(a) \dd y.
\end{split}
\end{equation}

\begin{proposition} \label{*lemintest}
  Let $b \ge 2$ be an integer, $0<\alpha<1$, and $W \in \WW^{\alpha,b}$. Then for $\LL^d$-almost every $\ttt \in \R^d$,
  \begin{equation*}
    \dimh(W_\ttt^{-1}(\{y\})) \ge 1-\alpha
  \end{equation*}
  for $\lambda_\ttt$-almost all $y \in \R$.
\end{proposition}

\begin{proof}
  If the $(1-\alpha)$-energy of $(\nu_{\ttt})_{y}$ is finite, i.e.,
  \begin{equation*}
    I_{1-\alpha}((\nu_{\ttt})_{y}) = \iint \|a-b\|^{-(1-\alpha)}\dd(\nu_{\ttt})_{y}(a)\dd(\nu_{\ttt})_{y}(b) < \infty,
  \end{equation*}
  then, by e.g.\ \cite[Theorem 4.13]{MR1102677}, the support of $(\nu_{\ttt})_{y}$ has Hausdorff dimension at least $1-\alpha$, that is, $\dimh(\graph(W_{\ttt}) \cap \proj_{2}^{-1}(\{y\})) \ge 1-\alpha$. Thus, it suffices to show that 
  \begin{equation} \label{*eqlemintest}
    \int\limits_{B(\ttt_0,\roo)}  \int I_{1-\beta}((\nu_{\ttt})_{y}) \dd y\dd\textbf{t} < \infty
  \end{equation}
  for all $\beta > \alpha$, $\ttt_0 \in \R^d$, and $\roo > 0$.

  By applying Fatou's lemma, \eqref{*disint}, and  \cref{eq:grslbcl6}, we obtain
  \begin{align*} %\label{eq:fatou-and-friends}
    \int\limits_{B(\ttt_0,\roo)}\int &I_{1-\beta}((\nu_{\ttt})_{y}) \dd y\dd\textbf{t} \\ 
    &= \int\limits_{B(\ttt_0,\roo)}  \iiint \frac{1}{\|a-b\|^{1-\beta}} \dd(\nu_{\textbf{t}})_{y}(a) \dd(\nu_{\textbf{t}})_{y}(b) \dd y\dd\textbf{t} \\ 
    &= \int\limits_{B(\ttt_0,\roo)}  \iiint \frac{h_{\ttt}(y)}{\|a-b\|^{1-\beta}} \dd(\mu_{\textbf{t}})_{y}(a) \dd(\nu_{\textbf{t}})_{y}(b) \dd y\dd\textbf{t} \\ 
    &\le \liminf_{r \downarrow 0} \frac{1}{2r} \int\limits_{B(\ttt_0,\roo)}  \iint \int\limits_{B^\top(y,r)} \frac{\mathrm{d}\mu_{\ttt}(a)}{\|a-b\|^{1-\beta}} \dd(\nu_{\textbf{t}})_{y}(b) \dd y\dd\ttt. 
  \end{align*}
  Hence, keeping in mind that $a\in B^\top(y,r)$ is equivalent to $y\in B(\proj_{2}(a),r)$, switching the order of integration, using \cref{eq:grslbcl5}, noting that $b\in B^{\top}(\proj_{2}(a),r)$ is equivalent to $a\in B^{\top}(\proj_{2}(b),r)$, and applying Fubini once more, we obtain
  \begin{equation}\label{eq:fatou-and-friendsb}
  \begin{split}
    \int\limits_{B(\ttt_0,\roo)}\int &I_{1-\beta}((\nu_{\ttt})_{y}) \dd y\dd\textbf{t} \\ 
    &\leq \liminf_{r \downarrow 0} \frac{1}{2r} \int\limits_{B(\ttt_0,\roo)} \iint_{\proj_{2}(a)-r}^{\proj_{2}(a)+r} \int   \frac{1}{\|a-b\|^{1-\beta}} \dd(\nu_{\textbf{t}})_{y}(b) \dd y \mathrm{d}\mu_{\ttt}(a) \dd\ttt \\ 
    &= \liminf_{r \downarrow 0} \frac{1}{2r} \int\limits_{B(\ttt_0,\roo)} \int \int\limits_{B^{\top}(\proj_{2}(a),r)}   \frac{1}{\|a-b\|^{1-\beta}}  \mathrm{d}\mu_{\ttt}(b) \mathrm{d}\mu_{\ttt}(a) \dd\ttt \\ 
    &= \liminf_{r \downarrow 0} \frac{1}{2r} \int\limits_{B(\ttt_0,\roo)}  \int \int\limits_{B^\top(\proj_{2}(b),r)} \frac{1}{\|a-b\|^{1-\beta}} \dd\mu_{\ttt}(a)\dd\mu_{\textbf{t}}(b)\dd\ttt.
  \end{split}
  \end{equation}
	Writing $a=(u,W_{\ttt}(u))$ and $b=(v,W_{\ttt}(v))$ with $u,v\in [0,1]$, \cref{eq:fatou-and-friendsb} and Fubini's theorem yield
  \begin{equation} \label{eq:inner-integral}
  \begin{split}
    \int\limits_{B(\ttt_0,\roo)}\int &I_{1-\beta}((\nu_{\ttt})_{y}) \dd y\dd\textbf{t} \\ 
	  &\leq \liminf_{r \downarrow 0} \frac{1}{2r} \int_{0}^{1}\int_{0}^{1} \int\limits_{B(\ttt_0,\roo)}\frac{\mathds{1}_{B^\top(W_{\ttt}(v),r)}(u,W_{\ttt}(u))}{\|(u,W_{\ttt}(u))-(v,W_{\ttt}(v))\|^{1-\beta}} \dd\ttt \dd u \dd v.
  \end{split}
  \end{equation}
  As the first step in estimating the inner integral in the upper bound of \cref{eq:inner-integral}, note that $\mathds{1}_{B^\top(W_{\ttt}(v),r)}(u,W_{\ttt}(u))=1$ if and only if $(u,W_{\ttt}(u))\in B^\top(W_{\ttt}(v),r)$, that is, 
  \begin{equation} \label{*W12*a}
    |W_{\ttt}(u)-W_{\ttt}(v)| < r.
  \end{equation}
  Recall that $W_{\ttt}(x) = W(x)+\langle \ttt, \Phi(x)\rangle$, where $\Phi = (W_0,W_1,\ldots,W_{d-1})$ is an $\alpha$-bi-H\"older Weierstrass embedding. Since $\|\cdot\|$ is the $\ell^\infty$ norm, we have
  \begin{equation}\label{*normest}
    \|(u,W_{\ttt}(u))-(v,W_{\ttt}(v))\| = \max\{ |u-v|, |W_{\ttt}(u)-W_{\ttt}(v)| \} \geq |u-v|
  \end{equation}
  and that there exists $i \in \{0,\ldots,d-1\}$ such that $\|\Phi(u)-\Phi(v)\| = |W_{i}(u)-W_{i}(v)|$. Without loss of generality, assume $i=0$. Define $\overline{\Phi} = (W_1,\ldots,W_{d-1})$ and 
for $\ttt \in \R^d$  write $\ttt = (t,\overline{\ttt}) \in \R \times \R^{d-1}$. Thus, \cref{*W12*a} is equivalent to
  \begin{equation*}
    |W(u)-W(v) + t\sgn(u,v)\|\Phi(u)-\Phi(v)\| + \langle \overline{\ttt}, \overline{\Phi}(u)-\overline{\Phi}(v) \rangle| < r,
  \end{equation*}
  where $\sgn(u,v) = (W_0(u)-W_0(v))/\|\Phi(u)-\Phi(v)\|
  = (W_0(u)-W_0(v))/|W_0(u)-W_0(v)|$. Define
  \begin{equation*}
    \tau(\overline{\ttt}) = -\frac{W(u)-W(v) + \langle \overline{\ttt}, \overline{\Phi}(u)-\overline{\Phi}(v) \rangle}{\sgn(u,v)\|\Phi(u)-\Phi(v)\|}.
  \end{equation*}
  Hence, \cref{*W12*a} holds if and only if
  \begin{equation} \label{*W15*a}
    t \in B\biggl( \tau(\overline{\ttt}), \frac{r}{\|\Phi(u)-\Phi(v)\|} \biggr).
  \end{equation}
  Having established these observations, we proceed to estimate the inner integral in the upper bound of \cref{eq:inner-integral}. Write $\ttt_0 = (t_{0},\overline{\ttt}_0) \in \R \times \R^{d-1}$. By \cref{*normest}, \cref{*W15*a}, and the bi-H\"older property of $\Phi$, we obtain
  \begin{align*}
    \int\limits_{B(\ttt_0,\roo)} &\frac{\mathds{1}_{B^\top(W_{\ttt}(v),r)}(u,W_{\ttt}(u))}{\|(u,W_{\ttt}(u))-(v,W_{\ttt}(v))\|^{1-\beta}} \dd\ttt \\ 
    &= \int\limits_{B(\overline{\ttt}_0,\roo)} \int\limits_{t_{0}-\roo}^{t_{0}+\roo} \frac{\mathds{1}_{B^\top(W_{(t,\overline{\ttt})}(v),r)}(u,W_{(t,\overline{\ttt})}(u))}{\|(u,W_{(t,\overline{\ttt})}(u))-(v,W_{(t,\overline{\ttt})}(v))\|^{1-\beta}} \dd t \dd\overline{\ttt}, \\ 
    &\leq \int\limits_{B(\overline{\ttt}_0,\roo)}\int\limits_{\R}\frac{\mathds{1}_{B^\top(W_{(t,\overline{\ttt})}(v),r)}(u,W_{(t,\overline{\ttt})}(u))}{\|(u,W_{(t,\overline{\ttt})}(u))-(v,W_{(t,\overline{\ttt})}(v))\|^{1-\beta}} \dd t \dd\overline{\ttt} \\ 
    &\leq \int\limits_{B(\overline{\ttt}_0,\roo)}\int\limits_{\tau(\overline{\ttt})-\frac{r}{\|\Phi(u)-\Phi(v)\|}}^{\tau(\overline{\ttt})+\frac{r}{\|\Phi(u)-\Phi(v)\|}}\frac{1}{|u-v|^{1-\beta}} \dd t \dd\overline{\ttt}. 
\end{align*}
  Hence with some constants $C'$ and $C$ not depending on $u$ and $v$,
  \begin{align*}	
    \int\limits_{B(\ttt_0,\roo)} &\frac{\mathds{1}_{B^\top(W_{\ttt}(v),r)}(u,W_{\ttt}(u))}{\|(u,W_{\ttt}(u))-(v,W_{\ttt}(v))\|^{1-\beta}} \dd\ttt \\ 
    &\leq \frac{1}{|u-v|^{1-\beta}}C'\roo^{d-1}\frac{r}{\|\Phi(u)-\Phi(v)\|} \\ 
    &\leq \frac{1}{|u-v|^{1-\beta}}C\roo^{d-1}\frac{r}{|u-v|^{\alpha}}=\frac{C\roo^{d-1} r}{|u-v|^{1-\beta+\alpha}}.
  \end{align*}
   Since $1-\beta+\alpha<1$, by substituting this estimate into \cref{eq:inner-integral}, we obtain
  \begin{equation*}
    \int\limits_{B(\ttt_0,\roo)}\int I_{1-\beta}((\nu_{\ttt})_{y}) \dd y\dd\textbf{t} \leq \int_{0}^{1}\int_{0}^{1}\frac{C\roo^{d-1}}{|u-v|^{1-\beta+\alpha}}\dd u\dd v < \infty.
  \end{equation*}
  This establishes \eqref{*eqlemintest} and completes the proof.
\end{proof}

The space
\begin{equation*}
  \SS = \{x\mapsto \langle\ttt, \Phi(x)\rangle : \ttt \in \R^d\}
\end{equation*}
is a $d$-dimensional subspace of $\WW^{\alpha,b}$ and serves as the probe space in our results. To establish that a property holds for $d$-prevalent functions, it suffices to demonstrate that for any $W\in\WW^{\alpha,b}$, the property holds for $W_{\ttt}$ for $\LL^d$-almost every $\ttt\in\R^d$.

\begin{proof}[Proof of \cref{thm:main2}]
  Let $b \geq 2$ an integer, $0<\alpha<1$, and let $\Phi=\Phi_{\mathcal{G}}^{\alpha,b}$ be the $\alpha$-bi-H\"older Weierstrass embedding given by \cref{prop:bi-holder}. By \cref{*lemintest} and \cref{*prabscb}, for every $W\in\WW^{\alpha,b}$ and $\LL^d$-almost every $\ttt \in \R^d$, the function $W_{\ttt} \in \mathcal{W}^{\alpha,b}$ defined by $W_{\ttt}(x)=W(x)+\langle\ttt,\Phi(x)\rangle$ satisfies
  \begin{equation*}
    \dimh(W_\ttt^{-1}(\{y\})) \ge 1-\alpha
  \end{equation*}
  for Lebesgue positively many $y \in \R$. The upper bound follows from \cref{eq:marstrand}.
\end{proof}

\section{Open questions}\label{secq}

In this section, we outline open problems arising from our work for future investigation.

\begin{enumerate}
  \item\label{q:q1} The primary open question is to strengthen \cref{thm:main2} by proving that for prevalent functions in $\mathcal{W}^{\alpha,b}$, with $b \geq 2$ and $0 < \alpha < 1$, the Hausdorff dimension of the level set $(W_g^{\alpha,b})^{-1}(\{y\})$ equals $1 - \alpha$ for Lebesgue almost all $y \in W_g^{\alpha,b}(\mathbb{S}^1)$. Additionally, could the Hausdorff dimension $\dimh((W_g^{\alpha,b})^{-1}(\{y\}))$ equal $1 - \alpha$ for all $y \in \mathrm{int}(W_g^{\alpha,b}(\mathbb{S}^1))$? To establish the Lebesgue almost all claim, it suffices to show that for prevalent functions $W_g^{\alpha,b} \in \mathcal{W}^{\alpha,b}$, the occupation measure has a density function 
which  is positive Lebesgue almost everywhere on $W_g^{\alpha,b}(\mathbb{S}^1)$. Since $\graph(W_g^{\alpha,b})$ is the attractor of an iterated function system defined by the maps
  \begin{equation*}
    \fii_j(x, y) = \Bigl( \frac{x + j}{b}, b^{-\alpha} y + g\Bigl( \frac{x + j}{b} \Bigr) \Bigr), \quad j \in \{0, \ldots, b-1\},
  \end{equation*}
  this could be achieved by analyzing the Lebesgue positively many level sets with absolutely continuous occupation measures through iterative application of these maps.

  \item\label{q:q2} Further investigation into the properties of the density function of the occupation measure would be valuable. For instance, extending \cref{*prabscbb} to $\frac12 \leq \alpha < 1$ would be of interest. Additionally, verifying 
  for the density function
  stronger properties than continuity  across certain ranges of $\alpha$ is a compelling challenge. Heuristically, there appears to be a dichotomy: when $\alpha$ is small, the prevalent Weierstrass function exhibits ``wilder'' behavior with a smaller H\"older exponent, yet its occupation measure's density function may possess more favorable properties.

  \item\label{q:q3} Instead of focusing solely on the occupation measure, one can consider slices and projections in various directions. For $\theta \in [0, 2\pi)$, let ${\mathrm {pr}}_\theta \colon \mathbb{R}^2 \to \mathbb{R}$ denote the orthogonal projection defined by
  \begin{equation*}
    {\mathrm {pr}}_\theta(x, y) = x \cos(\theta) + y \sin(\theta)
  \end{equation*}
  for all $(x, y) \in \mathbb{R}^2$. Note that $\proj_2 = {\mathrm {pr}}_{\frac{\pi}{2}}$. Given the measure $\mu_\ttt$ defined earlier in \cref{eq:lift-measure}, we consider the projections of $\mu_\ttt$ for $\theta \in [0, 2\pi)$, defined as
  \begin{equation*}
    \lambda_\ttt^\theta = ({\mathrm {pr}}_\theta)_* \mu_\ttt.
  \end{equation*}
  These measures $\lambda_\ttt^\theta$ are referred to as \emph{$\theta$-oblique occupation measures}. Note that the occupation measure is $\lambda_\ttt = \lambda_\ttt^{\frac{\pi}{2}}$. 

  By combining techniques from this paper and \cite{AnttilaBaranyKaenmaki2025}, one can show that for prevalent functions in $\mathcal{W}^{\alpha,b}$, the $\theta$-oblique occupation measure is absolutely continuous with respect to the Lebesgue measure for Lebesgue almost every $\theta$. A compelling open problem is to determine whether this property holds for all $\theta$, i.e., whether the $\theta$-oblique occupation measure is absolutely continuous for every $\theta$. Furthermore, one could explore whether questions analogous to \cref{q:q1} and \cref{q:q2} admit positive answers for all $\theta$ in this context. In particular, a stronger version of \cref{q:q1} could be investigated: can it be shown that for prevalent functions in $\mathcal{W}^{\alpha,b}$, the non-empty, or non-extremal slices in all non-vertical directions simultaneously have Hausdorff dimension $1 - \alpha$?

  \item \cref{prop:bi-holder} ensures the existence of $\alpha$-bi-H\"older Weierstrass embeddings only when the dimension $d$ is sufficiently large. 
  In \cref{*propaaabound}, a lower bound on $d$ is established as a function of $\alpha$. Is this bound sharp? What is the minimal dimension $d$ required for a given $\alpha$?
\end{enumerate}

\bibliographystyle{abbrv}
\bibliography{wlevBibliography}

@article {AlKa,
    AUTHOR = {Allaart, Pieter C. and Kawamura, Kiko},
     TITLE = {The {T}akagi function: a survey},
   JOURNAL = {Real Anal. Exchange},
  FJOURNAL = {Real Analysis Exchange},
    VOLUME = {37},
      YEAR = {2011/12},
    NUMBER = {1},
     PAGES = {1--54},
      ISSN = {0147-1937,1930-1219},
   MRCLASS = {26A27 (26A16 28A80)},
  MRNUMBER = {3016850},
MRREVIEWER = {Roland\ Girgensohn},
       URL = {http://projecteuclid.org/euclid.rae/1335806762},
}

@incollection {Lags,
    AUTHOR = {Lagarias, Jeffrey C.},
     TITLE = {The {T}akagi function and its properties},
 BOOKTITLE = {Functions in number theory and their probabilistic aspects},
    SERIES = {RIMS K\^oky\^uroku Bessatsu},
    VOLUME = {B34},
     PAGES = {153--189},
 PUBLISHER = {Res. Inst. Math. Sci. (RIMS), Kyoto},
      YEAR = {2012},
   MRCLASS = {26A27 (11A63)},
  MRNUMBER = {3014845},
MRREVIEWER = {Kiko\ Kawamura},
}

@article {GeHo,
AUTHOR = {Geman, Donald and Horowitz, Joseph},
TITLE = {Occupation densities},
JOURNAL = {Ann. Probab.},
FJOURNAL = {The Annals of Probability},
VOLUME = {8},
YEAR = {1980},
NUMBER = {1},
PAGES = {1--67},
ISSN = {0091-1798,2168-894X},
MRCLASS = {60J55 (26A27 60G15 60G17)},
MRNUMBER = {556414},
MRREVIEWER = {Simeon\ M.\ Berman},
URL =
{http://links.jstor.org/sici?sici=0091-1798(198002)8:1<1:OD>2.0.CO;2-M&origin=MSN},
}

@incollection {BZom,
AUTHOR = {Buczolich, Zolt\'an},
TITLE = {Occupation measure and level sets of the
{W}eierstrass-{C}ellerier function},
BOOKTITLE = {Recent developments in fractals and related fields},
SERIES = {Appl. Numer. Harmon. Anal.},
PAGES = {3--18},
PUBLISHER = {Birkh\"auser Boston, Boston, MA},
YEAR = {2010},
ISBN = {978-0-8176-4887-9},
MRCLASS = {28A80 (26A30)},
MRNUMBER = {2742982},
MRREVIEWER = {Pieter\ C.\ Allaart},
DOI = {10.1007/978-0-8176-4888-6\_1},
URL = {https://doi.org/10.1007/978-0-8176-4888-6_1},
}

@article {BZirr,
AUTHOR = {Buczolich, Z.},
TITLE = {Irregular 1-sets on the graphs of continuous functions},
JOURNAL = {Acta Math. Hungar.},
FJOURNAL = {Acta Mathematica Hungarica},
VOLUME = {121},
YEAR = {2008},
NUMBER = {4},
PAGES = {371--393},
ISSN = {0236-5294,1588-2632},
MRCLASS = {28A75 (26A24 26A27)},
MRNUMBER = {2461441},
MRREVIEWER = {Ville\ Suomala},
DOI = {10.1007/s10474-008-7220-9},
URL = {https://doi.org/10.1007/s10474-008-7220-9},
}

@article {Berom,
AUTHOR = {Bertoin, Jean},
TITLE = {Sur la mesure d'occupation d'une classe de fonctions
self-affines},
JOURNAL = {Japan J. Appl. Math.},
FJOURNAL = {Japan Journal of Applied Mathematics},
VOLUME = {5},
YEAR = {1988},
NUMBER = {3},
PAGES = {431--439},
ISSN = {0910-2043},
MRCLASS = {26A30 (11K55 60G17)},
MRNUMBER = {965873},
MRREVIEWER = {Michael\ Evans},
DOI = {10.1007/BF03167910},
URL = {https://doi.org/10.1007/BF03167910},
}

@article {Berhdlev,
AUTHOR = {Bertoin, Jean},
TITLE = {Hausdorff dimension of the level sets for self-affine
functions},
JOURNAL = {Japan J. Appl. Math.},
FJOURNAL = {Japan Journal of Applied Mathematics},
VOLUME = {7},
YEAR = {1990},
NUMBER = {2},
PAGES = {197--202},
ISSN = {0910-2043},
MRCLASS = {28A78 (60G17)},
MRNUMBER = {1057529},
MRREVIEWER = {K.\ E.\ Hirst},
DOI = {10.1007/BF03167841},
URL = {https://doi.org/10.1007/BF03167841},
}

@book{MattilaFou,
	author = {Mattila, P.},
	doi = {10.1017/CBO9781316227619},
	isbn = {978-1-107-10735-9},
	mrclass = {28-02 (28A15 28A78 28A80 42B10 60J65)},
	mrnumber = {3617376},
	mrreviewer = {Benjamin Steinhurst},
	pages = {xiv+440},
	publisher = {Cambridge University Press, Cambridge},
	series = {Cambridge Studies in Advanced Mathematics},
	title = {Fourier analysis and {H}ausdorff dimension},
	url = {https://doi.org/10.1017/CBO9781316227619},
	volume = {150},
	year = {2015}}

@article{AnttilaBaranyKaenmaki2025,
    AUTHOR = {Anttila, Roope and B\'ar\'any, Bal\'azs and K\"aenm\"aki,
              Antti},
     TITLE = {Level sets of prevalent {H}\"older functions},
   JOURNAL = {Proc. Amer. Math. Soc.},
  FJOURNAL = {Proceedings of the American Mathematical Society},
    VOLUME = {153},
      YEAR = {2025},
    NUMBER = {5},
     PAGES = {2023--2035},
      ISSN = {0002-9939,1088-6826},
   MRCLASS = {26A16 (26E15 28A50 28A78)},
  MRNUMBER = {4881392},
       DOI = {10.1090/proc/17045},
       URL = {https://doi.org/10.1090/proc/17045},
}

@article {Simmons2012,
    AUTHOR = {Simmons, D.},
     TITLE = {Conditional measures and conditional expectation; {Rohlin's} disintegration theorem},
   JOURNAL = {Disc. Cont. Dyn. Sys.},
    VOLUME = {32},
      YEAR = {2012},
    NUMBER = {7},
     PAGES = {2565--2582},
}

@book{BishopPeres2016,
    AUTHOR = {Bishop, C. J. and Peres, Y.},
     TITLE = {Fractals in Probability and Analysis},
    SERIES = {Universitext},
 PUBLISHER = {Cambridge: Cambridge University Press},
      YEAR = {2016},
}

@Book{ Mattila1995,
	address = "Cambridge",
	author = "Pertti Mattila",
	publisher = "Cambridge University Press",
	title = "Geometry of Sets and Measures in Euclidean Spaces: Fractals and Rectifiability",
	year = "1995"
}

@article {AnttilaBaranyKaenmaki2023,
    AUTHOR = {Anttila, Roope and Bárány, Balázs and K{\"a}enm{\"a}ki, Antti},
     TITLE = {Slices of the {T}akagi function},
   JOURNAL = {Ergodic Theory Dynam. Systems},
      YEAR = {2024},
    VOLUME = {44},
    NUMBER = {9},
     PAGES = {2361--2398},
}

@book {MR1102677,
    AUTHOR = {Falconer, Kenneth},
     TITLE = {Fractal geometry: Mathematical foundations and applications},
 PUBLISHER = {John Wiley \& Sons, Ltd., Chichester},
      YEAR = {2014},
   EDITION = {Third},
}

@article {Hochman2014,
    AUTHOR = {Hochman, M.},
     TITLE = {On self-similar sets with overlaps and inverse theorems for entropy},
   JOURNAL = {Ann. of Math.},
    VOLUME = {180},
    NUMBER = {2},
      YEAR = {2014},
     PAGES = {773--822},
}

@article {Shen2018,
    AUTHOR = {Shen, W.},
     TITLE = {{Hausdorff} dimension of the graphs of the classical {Weierstrass} functions},
   JOURNAL = {Math. Z.},
    VOLUME = {289},
      YEAR = {2018},
     PAGES = {223--266},
}

@article {RenShen2021,
    AUTHOR = {Ren, H. and Shen, W.},
     TITLE = {A dichotomy for the {Weierstrass}-type functions},
   JOURNAL = {Invent. Math.},
    VOLUME = {226},
      YEAR = {2021},
     PAGES = {1057--1100},
}

@InProceedings{10.1007/978-3-319-18660-3_5,
author="Bara{\'{n}}ski, Krzysztof",
editor="Bandt, Christoph
and Falconer, Kenneth
and Z{\"a}hle, Martina",
title="Dimension of the Graphs of the {W}eierstrass-Type Functions",
booktitle="Fractal Geometry and Stochastics V",
year="2015",
publisher="Springer International Publishing",
address="Cham",
pages="77--91",
isbn="978-3-319-18660-3"
}

@book {Morters_Peres_2010,
    AUTHOR = {M\"orters, Peter and Peres, Yuval},
     TITLE = {Brownian motion},
    SERIES = {Cambridge Series in Statistical and Probabilistic Mathematics},
    VOLUME = {30},
      NOTE = {With an appendix by Oded Schramm and Wendelin Werner},
 PUBLISHER = {Cambridge University Press, Cambridge},
      YEAR = {2010},
     PAGES = {xii+403},
      ISBN = {978-0-521-76018-8},
   MRCLASS = {60J65 (28A78 60H05 60J45 60J55 60J67)},
  MRNUMBER = {2604525},
MRREVIEWER = {Ren\'e\ L.\ Schilling},
       URL = {https://doi.org/10.1017/CBO9780511750489},
}

\end{document}